\documentclass[aop,preprint]{imsart}

\RequirePackage{amsthm,amsmath,amsfonts,amssymb,mathtools}
\RequirePackage[numbers,sort&compress]{natbib}
\RequirePackage[colorlinks,citecolor=blue,urlcolor=blue]{hyperref}

\startlocaldefs
\numberwithin{equation}{section}

\theoremstyle{plain}
\newtheorem{theorem}{Theorem}[section]
\newtheorem{lemma}[theorem]{Lemma}
\newtheorem{corollary}[theorem]{Corollary}
\newtheorem{proposition}[theorem]{Proposition}

\theoremstyle{definition}
\newtheorem{remark}[theorem]{Remark}
\newtheorem*{remark*}{Remark}

\newcommand{\R}{\mathbb{R}}

\newcommand{\E}{\mathbb{E}}
\newcommand{\Pp}{\mathbb{P}}
\newcommand{\B}{\mathcal{B}}

\newcommand{\law}{\mathcal{L}}
\newcommand{\ind}{\mathbf{1}}
\newcommand{\wh}{\widehat}
\newcommand{\dd}{\,\mathrm{d}}
\newcommand{\Var}{\operatorname{Var}}
\newcommand{\Cov}{\operatorname{Cov}}
\newcommand{\iid}{\mathrm{i.i.d.}}
\newcommand{\Normal}{\mathcal{N}}

\endlocaldefs

\begin{document}

\begin{frontmatter}

\title{Gaussian rigidity for infinite exchangeable sequences}
\runtitle{Gaussian rigidity for exchangeable sequences}
\runauthor{Y. Zheng and Q. Zhou}

\begin{aug}
\author[A]{\fnms{Yushu}~\snm{Zheng}\ead[label=e1]{yszheng666@gmail.com}}
\author[B]{\fnms{Qi}~\snm{Zhou}\ead[label=e2]{zhouqi@suda.edu.cn}}
\address[A]{Academy of Mathematics and Systems Science, Chinese Academy of Sciences\printead[presep={,\ }]{e1}}
\address[B]{School of Mathematical Sciences, Soochow University\printead[presep={,\ }]{e2}}
\end{aug}

\begin{abstract}
We prove a Gaussian rigidity theorem for infinite exchangeable sequences of
real-valued random variables: the joint Gaussianity of a single pair of
entries already forces the entire sequence to be a Gaussian process.  This
settles a conjecture raised by Newman (2026).  The main analytic
ingredient in the proof is Hardy's uncertainty principle.  We also obtain a
finite-dimensional vector-valued extension.
\end{abstract}



\end{frontmatter}

\section{Introduction}

Gaussian processes are defined by their finite-dimensional distributions, so
checking Gaussianity is usually a global task.  A natural question is whether,
under suitable structural assumptions, this global property can be detected
from much smaller pieces of the process.  The present paper proves such a
result for infinite exchangeable sequences: the joint Gaussianity of one pair of distinct coordinates already forces the whole sequence to be a Gaussian process.

Recall that an infinite sequence \((X_i)_{i\ge1}\) of \(\R^d\)-valued random
vectors is exchangeable if its law is invariant under finite permutations
of the indices. By the de Finetti--Hewitt--Savage representation theorem,
every infinite exchangeable sequence is conditionally i.i.d. given a random
probability measure, usually called its directing measure; see de Finetti~\cite{deFinetti1937}, Hewitt and Savage~\cite[Theorem~7.4]{HewittSavage1955}, and Kallenberg~\cite[Theorem~1.1]{Kallenberg2005}.

Newman~\cite{Newman2026} recently formulated the following problem for real-valued exchangeable sequences: if an infinite exchangeable sequence has a single jointly Gaussian pair $(X_1,X_2)$, must the whole sequence be a Gaussian process? This may be viewed as a two-point test for Gaussianity within the class of infinite exchangeable sequences. Newman proved a four-point version of this result, and also proved the two-point statement under the additional assumption that the directing measure is almost surely infinitely divisible.

\subsection{Main results}

We prove Newman's two-point conjecture in full, with no infinite-divisibility
assumption. The statement is as follows.

\begin{theorem}\label{thm:main}
Let $(X_i)_{i\ge1}$ be an infinite exchangeable sequence of real-valued random variables. If $(X_1,X_2)$ is jointly Gaussian, then $(X_i)_{i\ge1}$ is a Gaussian process.
\end{theorem}

A finite-dimensional vector-valued extension is proved in
Theorem~\ref{thm:vector-extension}.

\begin{remark*}
\begin{longlist}[(ii)]
\item It is equivalent to assume that there exist distinct indices
    $i\ne j$ such that $(X_i,X_j)$ is jointly Gaussian.  Indeed, by
    exchangeability, $(X_i,X_j)$ has the same distribution as $(X_1,X_2)$
    whenever $i\ne j$.

\item The corresponding one-point assumption is not sufficient.  For
    example, let $G\sim \Normal(0,1)$, let $R=|G|$, and let
    $\varepsilon_1,\varepsilon_2,\ldots$ be i.i.d. Rademacher random
    variables, independent of $R$.  Set
    \[
        X_i=R\varepsilon_i,\qquad i\ge1.
    \]
    Then $(X_i)_{i\ge1}$ is exchangeable, and each $X_i$ has the standard
    Gaussian distribution.  However, the sequence is not a Gaussian process:
    indeed, $\operatorname{Cov}(X_1,X_2)=0$, while $X_1$ and $X_2$ are not
    independent since $|X_1|=|X_2|$ almost surely.
\end{longlist}
\end{remark*}

\subsection{Proof outline}

We now describe the main ideas of the proof. The starting
point is the standard representation of exchangeable Gaussian sequences,
recorded in Proposition~\ref{prop:exchangeable-gaussian-representation}:
such a sequence is, in distribution, a common Gaussian shift plus independent
Gaussian noises. Thus subtracting two coordinates removes the common shift.
In de Finetti language, this suggests the following key identity for the
directing measure \(\eta\) (see Section~\ref{sec:definetti-identities} for
precise definition) of an exchangeable Gaussian sequence:
\[
        \eta*\widetilde\eta
        \quad\text{is a fixed centered Gaussian law a.s.},
\]
where \(\widetilde\eta(B)=\eta(-B)\).  Proposition~\ref{prop:gaussian-equivalences} shows that Theorem~\ref{thm:main} is reduced precisely to proving this identity.

The proof of this identity is where Hardy's uncertainty principle (see Theorem~\ref{thm:hardy}) enters. Hardy's theorem says that a function and its Fourier transform cannot both have matching Gaussian decay at the critical scale unless the function itself is a Gaussian. After localizing $\eta*\widetilde\eta$ on events in $\sigma(\eta)$, one obtains finite measures dominated by the Gaussian law of $X_1-X_2$. The equality \( \widehat{\eta*\widetilde\eta}(t)=|\widehat\eta(t)|^2 \) gives the matching Gaussian bound on the Fourier side. Hardy's theorem then forces each localization to be a scalar multiple of the same Gaussian law, and varying the event yields the desired almost sure identity for $\eta*\widetilde\eta$. The details are given in Section~\ref{sec:proof-difference-law}.

The vector-valued extension, Theorem~\ref{thm:vector-extension}, follows by applying Theorem~\ref{thm:main} to one-dimensional projections.

\subsection{Organization}

The rest of the paper is organized as follows.  Section~\ref{sec:preliminaries} records the preliminary ingredients: Hardy's uncertainty principle and de Finetti identities.  Section~\ref{sec:scalar-case} proves Theorem~\ref{thm:main}.  Section~\ref{sec:vector-extension} proves the vector-valued extension.

\section{Preliminaries}\label{sec:preliminaries}

In this section we record the preliminary ingredients used in the proof: Hardy's uncertainty principle and the elementary identities coming from de Finetti's representation.

\subsection{Hardy's uncertainty principle}

We use the Fourier transform convention
\[
        \widehat{f}(t)=\int_{\R} e^{itx} f(x)\dd x,
        \qquad
        \widehat{\nu}(t)=\int_{\R} e^{itx}\nu(\dd x).
\]
With this Fourier convention,
\[
        \int_{\R} e^{itx} e^{-x^2/(2a)}\dd x
        =\sqrt{2\pi a}\,e^{-a t^2/2}.
\]
Thus the Gaussian decay profiles \(e^{-x^2/(2a)}\) and \(e^{-a t^2/2}\) are Fourier-matched.

Hardy's uncertainty principle makes precise the idea that a function and its Fourier transform cannot both be too localized.  In the Gaussian scale, the product \(\alpha\beta=1/4\) is the borderline for decay profiles \(e^{-\alpha x^2}\) and \(e^{-\beta t^2}\); at this borderline, Hardy's theorem says that the matched Gaussian profile above is essentially the only nonzero possibility.

\begin{theorem}[Hardy's uncertainty principle]\label{thm:hardy}
Let $a>0$, and let $f$ be a measurable function such that
\[
        f(x)=O\!\left(e^{-x^2/(2a)}\right),\qquad
        \widehat f(t)=O\!\left(e^{-a t^2/2}\right)
\]
as $|x|\to\infty$ and $|t|\to\infty$, respectively.
Then there is a constant $\kappa\in\mathbb{C}$ such that
\[
        f(x)=\kappa e^{-x^2/(2a)}.
\]
\end{theorem}

\begin{remark}
This formulation is a simple rescaling of Hardy~\cite[Theorem~2]{Hardy1933}, where the theorem is stated in the normalization corresponding to \(a=1\); see also Fern\'andez-Bertolin and Vega~\cite[Remark~1]{FernandezBertolinVega2025}.
\end{remark}

We shall only need the following finite-measure corollary.

\begin{corollary}\label{cor:eventwise-hardy}
Let $a>0$, let $\gamma_a$ denote the centered Gaussian law $\Normal(0,a)$ on $\R$, and let $\nu$ be a finite positive Borel measure on $\R$.  Assume that
\begin{equation}\label{eq:measure-domination}
        0\le \nu\le \gamma_a
\end{equation}
as measures, and that
\begin{equation}\label{eq:fourier-domination}
        0\le \widehat\nu(t)\le e^{-a t^2/2}
        \qquad\text{for every }t\in\R.
\end{equation}
Then
\[
        \nu=\nu(\R)\,\gamma_a.
\]
\end{corollary}

\begin{proof}
By \eqref{eq:measure-domination}, the measure $\nu$ is absolutely continuous with respect to Lebesgue measure.  Write $\nu(\dd x)=h(x)\dd x$.  Since $\gamma_a$ has density
\[
        \varphi_a(x)=\frac{1}{\sqrt{2\pi a}}\exp\left\{-\frac{x^2}{2a}\right\},
\]
we may choose a version of \(h\), possibly after modifying it on a Lebesgue null set, such that
\[
        0\le h(x)\le \varphi_a(x)
        \quad\text{for every }x\in\R.
\]
In particular, \(h(x)=O(e^{-x^2/(2a)})\) as \(|x|\to\infty\).  Moreover, \(\widehat h(t)=\widehat\nu(t)\).
By \eqref{eq:fourier-domination},
\[
        |\widehat h(t)|\le e^{-a t^2/2}.
\]

Hardy's uncertainty principle, Theorem~\ref{thm:hardy}, therefore implies that
\[
        h(x)=\kappa e^{-x^2/(2a)}
\]
for some constant $\kappa\ge0$.  Since $h\ge0$ and $\int h(x)\dd x=\nu(\R)$, the constant is
\[
        \kappa=\frac{\nu(\R)}{\sqrt{2\pi a}}.
\]
Thus $\nu=\nu(\R)\gamma_a$.
\end{proof}

\subsection{De Finetti identities}\label{sec:definetti-identities}

By de Finetti's theorem, if \((X_i)_{i\ge1}\) is an exchangeable sequence of \(\R^d\)-valued random vectors, then there exists a random probability measure \(\eta\) on \(\R^d\) such that, conditionally on \(\eta\), the variables \(X_1,X_2,\ldots\) are i.i.d. with common law \(\eta\).  This random probability measure \(\eta\) is called the de Finetti directing measure of the sequence.

We shall use the following elementary consequences of the de Finetti representation.

\begin{lemma}\label{lem:definetti-identities}
Let \((X_i)_{i\ge1}\) be an exchangeable sequence of \(\R^d\)-valued random vectors with de Finetti directing measure \(\eta\).  For every Borel function \(\phi:\R^d\to\mathbb C\) such that \(\E|\phi(X_1)|<\infty\),
\[
        \E[\phi(X_1)\mid\eta]
        =\int_{\R^d}\phi(x)\eta(\dd x)
        \qquad\text{a.s.}
\]
In particular,
\[
        \E\int_{\R^d}\phi(x)\eta(\dd x)=\E\phi(X_1)
\]
Moreover, for Borel functions \(\phi,\psi:\R^d\to\mathbb C\) such that \(\E|\phi(X_1)\psi(X_2)|<\infty\),
\[
        \E\bigl[\phi(X_1)\psi(X_2)\mid\eta\bigr]
        =\int_{\R^d}\phi(x)\eta(\dd x)\int_{\R^d}\psi(y)\eta(\dd y)
        \qquad\text{a.s.}
\]
\end{lemma}

\begin{proof}
Conditionally on \(\eta\), the variables \(X_1,X_2,\ldots\) are i.i.d. with common law \(\eta\).  Taking conditional expectations gives the displayed identities.
\end{proof}

We shall often use the following specialization.  Writing
\[
        \widehat\eta(u):=\int_{\R^d}e^{iu^{\mathsf T}x}\eta(\dd x),
\]
Lemma~\ref{lem:definetti-identities}, applied with
\(\phi(x)=e^{iu^{\mathsf T}x}\) and \(\psi(y)=e^{-iu^{\mathsf T}y}\), gives
\begin{equation}\label{eq:definetti-fourier-difference}
        \E\bigl[e^{iu^{\mathsf T}(X_1-X_2)}\mid\eta\bigr]
        =|\widehat\eta(u)|^2,
        \qquad u\in\R^d.
\end{equation}

\section{Proof of Theorem~\ref{thm:main}}\label{sec:scalar-case}

In this section we prove Theorem~\ref{thm:main}.  Section~\ref{sec:scalar-second-moment} records basic second-moment facts, Section~\ref{sec:scalar-reduction} reduces the theorem to an identity for the difference
law, and Section~\ref{sec:proof-difference-law} proves this identity.

Throughout this section, \((X_i)_{i\ge1}\) denotes an infinite real-valued exchangeable process with de Finetti directing measure \(\eta\), and we write
\[
        m:=\E X_1,\qquad
        \sigma^2:=\Var(X_1),\qquad
        c:=\Cov(X_1,X_2).
\]

\subsection{Second-moment facts}\label{sec:scalar-second-moment}

We first record two elementary properties of exchangeable processes with finite second moments.

\begin{lemma}\label{lem:scalar-second-moment}
Suppose that \(\E X_1^2<\infty\).  Then \(M_\eta:=\int x\eta(\dd x)\) is well-defined a.s. and \(M_\eta\in L^2\).  Moreover,
\[
        0\le c\le\sigma^2.
\]
\end{lemma}

\begin{proof}
By Lemma~\ref{lem:definetti-identities}, applied with \(\phi(x)=x^2\),
\begin{equation}\label{eq:scalar-eta-second-moment}
        \E\int x^2\eta(\dd x)=\E X_1^2<\infty.
\end{equation}
Thus \(\int x^2\eta(\dd x)<\infty\) a.s., so \(M_\eta\) is well-defined a.s.  By Jensen's inequality,
\[
        M_\eta^2\le \int x^2\eta(\dd x),
\]
and hence \(M_\eta\in L^2\).

Lemma~\ref{lem:definetti-identities}, applied with \(\phi(x)=x\) and with \(\phi(x)=\psi(x)=x\), gives
\[
        \E[X_i\mid\eta]=M_\eta,\qquad
        \E[X_1X_2\mid\eta]=M_\eta^2.
\]
Therefore
\[
        c=\Cov(X_1,X_2)=\E(M_\eta^2)-\E(M_\eta)^2=\Var(M_\eta)\ge0.
\]
By the law of total variance,
\[
        \sigma^2=\E[\Var(X_1\mid\eta)]+\Var(M_\eta)\ge c.
\]
This proves the claim.
\end{proof}

\subsection{Reduction to the difference law}\label{sec:scalar-reduction}

The goal of this subsection is to reduce Theorem~\ref{thm:main} to identifying the difference law of the directing measure, as stated in \eqref{eq:target-difference-law}.

We first record the standard representation of Gaussian exchangeable sequences.
This representation goes back to
Aldous~\cite[p.~8]{Aldous1985}; see also
Newman~\cite[Equation~(1)]{Newman2026}.

\begin{proposition}\label{prop:exchangeable-gaussian-representation}
Assume that \((X_i)_{i\ge1}\) is a Gaussian exchangeable process.
The law of \((X_i)_{i\ge1}\) has the following representation in distribution. Let \(M,\xi_1,\xi_2,\ldots\) be independent real-valued random variables such that
\[
        M\sim\Normal(m,c),\qquad
        \xi_1,\xi_2,\ldots \text{ are i.i.d. } \Normal(0,\sigma^2-c),
\]
(the nonnegativity of \(c\) and \(\sigma^2-c\) follows from
Lemma~\ref{lem:scalar-second-moment}), and set \(Y_i=M+\xi_i\).  Then
\[
        (X_i)_{i\ge1}\stackrel d=(Y_i)_{i\ge1}.
\]
\end{proposition}

\begin{proof}
By exchangeability, for each \(n\), the covariance matrix of \((X_1,\ldots,X_n)\) is
\[
        \Sigma_n:=
        \begin{pmatrix}
        \sigma^2 & c & \cdots & c \\
        c & \sigma^2 & \cdots & c \\
        \vdots & \vdots & \ddots & \vdots \\
        c & c & \cdots & \sigma^2
        \end{pmatrix}.
\]

Let \(M,\xi_1,\xi_2,\ldots\) be independent random variables with the laws specified above, and set \(Y_i=M+\xi_i\).  For each \(n\), both \((X_1,\ldots,X_n)\) and \((Y_1,\ldots,Y_n)\) are Gaussian vectors with mean \(m\mathbf{1}_n\) and covariance matrix \(\Sigma_n\).
Hence the two sequences have the same finite-dimensional distributions.
\end{proof}

In the standard representation, the common random shift \(M\) disappears when two coordinates are subtracted, so the difference law is governed only by the noise variables.  Motivated by this observation, we next reduce Gaussianity to an identity for the difference law \(\eta*\widetilde\eta\).

\begin{proposition}\label{prop:gaussian-equivalences}
Suppose that \(X_1\sim\Normal(m,\sigma^2)\).
Set \(a:=2(\sigma^2-c)\), which is nonnegative by Lemma~\ref{lem:scalar-second-moment}, and let \(\gamma_a\) denote the centered Gaussian law \(\Normal(0,a)\) on \(\R\).  Then the following are equivalent.
\begin{enumerate}
\renewcommand{\labelenumi}{(\roman{enumi})}
\renewcommand{\theenumi}{\roman{enumi}}
\item\label{prop:gaussian-equivalences:i} \((X_i)_{i\ge1}\) is a Gaussian process.
\item\label{prop:gaussian-equivalences:ii} \(\eta*\widetilde\eta=\gamma_a\) a.s., where \(\widetilde\eta(B)=\eta(-B)\).
\item\label{prop:gaussian-equivalences:iii} \(M_\eta\sim\Normal(m,c)\) and \(\eta\) is the Gaussian law \(\Normal(M_\eta,\sigma^2-c)\) a.s.
\end{enumerate}
\end{proposition}

\begin{remark*}
Conditionally on \(\eta\), the variables \(X_1\) and \(X_2\) are i.i.d. with common law \(\eta\).  Hence \(\eta*\widetilde\eta\) is the conditional law of \(X_1-X_2\) given \(\eta\).
\end{remark*}

\begin{proof}
We prove the cycle \textnormal{(\ref{prop:gaussian-equivalences:i})}
\(\Rightarrow\)\textnormal{(\ref{prop:gaussian-equivalences:ii})}
\(\Rightarrow\)\textnormal{(\ref{prop:gaussian-equivalences:iii})}
\(\Rightarrow\)\textnormal{(\ref{prop:gaussian-equivalences:i})}.

Assume \textnormal{(\ref{prop:gaussian-equivalences:i})}.  By Proposition~\ref{prop:exchangeable-gaussian-representation}, the sequence admits the standard representation
\[
        (X_i)_{i\ge1}\stackrel d=(M+\xi_i)_{i\ge1},
\]
where \(M\sim\Normal(m,c)\), the variables \(\xi_i\) are i.i.d. \(\Normal(0,\sigma^2-c)\), and \(M\) is independent of \((\xi_i)_{i\ge1}\).  Conditionally on \(M\), the represented sequence is i.i.d. with common law \(\Normal(M,\sigma^2-c)\).  Thus, by uniqueness of the de Finetti directing measure (see, for instance, Kallenberg~\cite[Theorem~1.1]{Kallenberg2005}),
\[
        \eta\stackrel d=\Normal(M,\sigma^2-c).
\]
For this represented sequence, the difference law \(\eta*\widetilde\eta\) has the
same distribution as
\[
        \law\bigl((M+\xi_1)-(M+\xi_2)\mid M\bigr)
        =
        \law(\xi_1-\xi_2)
        =
        \Normal(0,2(\sigma^2-c))
        =
        \gamma_a,
\]
which is deterministic. Hence \(\eta*\widetilde\eta=\gamma_a\) a.s.
This proves \textnormal{(\ref{prop:gaussian-equivalences:ii})}.

Assume \textnormal{(\ref{prop:gaussian-equivalences:ii})}.  Fix a value of \(\eta\) for which \(\eta*\widetilde\eta=\gamma_a\).  If \(U,V\) are independent with common law \(\eta\), then \(U-V\sim\Normal(0,a)\).  By Cram\'er's decomposition theorem for normal distributions, allowing degenerate normal laws~\cite{Cramer1936}, both \(U\) and \(-V\) are Gaussian.  Since \(U\) and \(V\) have the same law and \(U-V\) has variance \(a\), this Gaussian law has variance \(a/2=\sigma^2-c\).  Hence \(\eta=\Normal(M_\eta,\sigma^2-c)\) a.s.  Using \(X_1\sim\Normal(m,\sigma^2)\), we get, for every \(t\in\R\),
\[
        \exp\left\{imt-\frac{\sigma^2t^2}{2}\right\}
        =\E e^{itX_1}
        =\E\left[\exp\left\{itM_\eta-\frac{(\sigma^2-c)t^2}{2}\right\}\right].
\]
Thus
\[
        \E e^{itM_\eta}
        =\exp\left\{imt-\frac{ct^2}{2}\right\},
        \qquad\text{so}\qquad
        M_\eta\sim\Normal(m,c).
\]
This proves \textnormal{(\ref{prop:gaussian-equivalences:iii})}.

Finally, assume \textnormal{(\ref{prop:gaussian-equivalences:iii})}.  By the conditional i.i.d. structure given \(\eta\), the finite-dimensional distributions of \((X_i)_{i\ge1}\) are those of \((M+\xi_i)_{i\ge1}\), where \(M\sim\Normal(m,c)\), the variables \(\xi_i\) are i.i.d. \(\Normal(0,\sigma^2-c)\), and \(M\) is independent of \((\xi_i)_{i\ge1}\).  Thus \((X_i)_{i\ge1}\) is a Gaussian process.
\end{proof}

We now turn to Theorem~\ref{thm:main}. Set \(a:=2(\sigma^2-c)\), and let \(\gamma_a\) denote the centered Gaussian law \(\Normal(0,a)\) on \(\R\).  By Proposition~\ref{prop:gaussian-equivalences}, it remains to prove
\begin{equation}\label{eq:target-difference-law}
        \eta*\widetilde\eta=\gamma_a
        \qquad\text{a.s.}
\end{equation}

\subsection{Proof of \textnormal{(\ref{eq:target-difference-law})}}\label{sec:proof-difference-law}

Since \((X_1,X_2)\) is bivariate Gaussian, \(X_1-X_2\sim\gamma_a\).  If \(a=0\), then \(X_1-X_2=0\) a.s.  Since the conditional law of \(X_1-X_2\) given \(\eta\) is \(\eta*\widetilde\eta\), it follows that \(\eta*\widetilde\eta=\delta_0=\gamma_0\) a.s.  Hence we may assume \(a>0\).

We localize the conditional difference law on events in \(\sigma(\eta)\).  For an arbitrary event \(E\in\sigma(\eta)\), define a finite positive measure \(\nu_E\) on \(\R\) by
\begin{equation}\label{eq:def-nu-E}
        \nu_E(B):=\E\bigl[\ind_E\cdot(\eta*\widetilde\eta)(B)\bigr]
        =\Pp\bigl(E\cap\{X_1-X_2\in B\}\bigr),
        \qquad B\in\B(\R).
\end{equation}

\begin{lemma}\label{lem:nuE-bounds}
With the above notation,
\[
        0\le \nu_E\le \gamma_a
\]
as measures, and
\[
        0\le \widehat{\nu_E}(t)\le e^{-a t^2/2},
        \qquad t\in\R.
\]
Consequently,
\begin{equation}\label{eq:eventwise-difference-law}
        \nu_E=\Pp(E)\gamma_a.
\end{equation}
\end{lemma}

\begin{proof}
For every Borel set $B$,
\[
        0\le \nu_E(B)
        \le \Pp(X_1-X_2\in B)
        =\gamma_a(B).
\]
This proves the measure domination.

For the Fourier bound,
\[
        \widehat{\nu_E}(t)
        =\int e^{itx}\nu_E(\dd x)
        \overset{\eqref{eq:def-nu-E}}=
        \E\bigl[\ind_E\cdot e^{it(X_1-X_2)}\bigr]
        \overset{\eqref{eq:definetti-fourier-difference}}=
        \E\bigl[\ind_E\cdot|\widehat\eta(t)|^2\bigr].
\]
The right-hand side is nonnegative and real.  Moreover,
\begin{equation*}
        0\le \widehat{\nu_E}(t)
        \le \E |\widehat\eta(t)|^2
        =\E e^{it(X_1-X_2)}
        =e^{-a t^2/2}.
\end{equation*}
The conclusion \eqref{eq:eventwise-difference-law} follows from Corollary~\ref{cor:eventwise-hardy}.
\end{proof}

By \eqref{eq:def-nu-E} and \eqref{eq:eventwise-difference-law}, for all $E\in\sigma(\eta)$ and $B\in\B(\R)$,
\[
        \E\bigl[\ind_E\cdot(\eta*\widetilde\eta)(B)\bigr]
        =\Pp(E)\gamma_a(B).
\]
Thus, for each fixed $B\in\B(\R)$,
\[
        (\eta*\widetilde\eta)(B)=\gamma_a(B)
        \qquad\text{a.s.}
\]
Applying this to a countable determining class, for instance the intervals $(-\infty,q]$ with $q\in\mathbb{Q}$, proves $\eta*\widetilde\eta=\gamma_a$ a.s.  This proves \eqref{eq:target-difference-law}, and hence Theorem~\ref{thm:main}.

\section{Vector-valued extension}\label{sec:vector-extension}

We now prove the finite-dimensional vector-valued extension of
Theorem~\ref{thm:main}. The proof applies the scalar result to
one-dimensional projections and then recovers the directing measure from
its characteristic function.

\begin{theorem}\label{thm:vector-extension}
Let $d \ge 1$, and let $(X_i)_{i\ge1}$ be an infinite exchangeable sequence
of $\mathbb{R}^d$-valued random vectors. Suppose that $(X_1,X_2)$ is Gaussian
as an $\mathbb{R}^{2d}$-valued random vector. Then $(X_i)_{i\ge1}$ is a
Gaussian process.
\end{theorem}

We start with the following vector-valued analogue of
Lemma~\ref{lem:scalar-second-moment}.

\begin{lemma}\label{lem:vector-second-moment}
Let \((X_i)_{i\ge1}\) be an exchangeable sequence of \(\R^d\)-valued random
vectors with de Finetti directing measure \(\eta\).  Suppose that
\(\E |X_1|^2<\infty\).  Then
\[
        M_\eta:=\int_{\R^d}x\eta(\dd x)
\]
is well-defined a.s. and belongs to \(L^2\).  Moreover, with
\(\Sigma:=\Cov(X_1)\) and \(C:=\Cov(X_1,X_2)\) (which is symmetric by exchangeability), we have
\[0\preceq C\preceq\Sigma.\]
Here and below, \(A\preceq B\) means that \(B-A\) is positive semidefinite.
\end{lemma}

\begin{proof}
The first assertion follows by the argument of
Lemma~\ref{lem:scalar-second-moment}, with \(|x|^2\) in place of \(x^2\).

For \(u\in\R^d\), the scalar
exchangeable sequence \((u^{\mathsf T}X_i)_{i\ge1}\) has variance
\(u^{\mathsf T}\Sigma u\) and two-point covariance \(u^{\mathsf T}Cu\).
Lemma~\ref{lem:scalar-second-moment} therefore gives
\[
        0\le u^{\mathsf T}Cu\le u^{\mathsf T}\Sigma u,
        \qquad u\in\R^d.
\]
Thus \(0\preceq C\preceq\Sigma\).
\end{proof}

\begin{proof}[Proof of Theorem~\ref{thm:vector-extension}]
Set \(m:=\E X_1\), \(\Sigma:=\Cov(X_1)\), and \(C:=\Cov(X_1,X_2)\).
Let \(\eta\) be the de Finetti directing measure on \(\R^d\).  Since \(X_1\)
is Gaussian, \(\E |X_1|^2<\infty\).  By Lemma~\ref{lem:vector-second-moment},
\(M_\eta:=\int_{\R^d}x\eta(\dd x)\) is well-defined a.s. and
\(0\preceq C\preceq\Sigma\).

Let \(D\subset\R^d\) be a countable dense set.  For each \(u\in D\), define
the scalar exchangeable sequence
\[
        Y_i^{(u)}:=u^{\mathsf T}X_i,\qquad i\ge1.
\]
Since \((X_1,X_2)\) is Gaussian as an \(\R^{2d}\)-valued random vector,
\((Y_1^{(u)},Y_2^{(u)})\) is bivariate Gaussian.  The de Finetti directing
measure of \((Y_i^{(u)})_{i\ge1}\) is the projected measure
\[
        \eta_u:=\eta\circ(x\mapsto u^{\mathsf T}x)^{-1}.
\]
Its mean is \(u^{\mathsf T}M_\eta\), while
\[
        \Var(Y_1^{(u)})=u^{\mathsf T}\Sigma u,
        \qquad
        \Cov(Y_1^{(u)},Y_2^{(u)})=u^{\mathsf T}Cu.
\]
By Theorem~\ref{thm:main}, the scalar sequence \((Y_i^{(u)})_{i\ge1}\) is a
Gaussian process.  Applying
Proposition~\ref{prop:gaussian-equivalences}\textnormal{(\ref{prop:gaussian-equivalences:iii})}
to this projected sequence gives
\[
        \eta_u
        =
        \Normal\bigl(u^{\mathsf T}M_\eta,
        u^{\mathsf T}(\Sigma-C)u\bigr)
        \qquad\text{a.s.}
\]
Since \(D\) is countable, almost surely, for every \(u\in D\),
\[
        \int_{\R^d}e^{iu^{\mathsf T}x}\eta(\dd x)
        =
        \exp\left\{
        iu^{\mathsf T}M_\eta
        -\frac12 u^{\mathsf T}(\Sigma-C)u
        \right\}.
\]
For each fixed value of \(\eta\) and \(M_\eta\), both sides are continuous
functions of \(u\in\R^d\).  Since \(D\) is dense, the identity extends to every
\(u\in\R^d\).  By uniqueness of characteristic functions on \(\R^d\),
\begin{equation}\label{eq:eta-vector-gaussian}
        \eta=\Normal(M_\eta,\Sigma-C)
        \qquad\text{a.s.}
\end{equation}

It remains to identify the law of \(M_\eta\).  Since \(X_1\sim\Normal(m,\Sigma)\)
and, conditionally on \(\eta\), \(X_1\) has law \(\Normal(M_\eta,\Sigma-C)\),
we have, for every \(u\in\R^d\),
\begin{align*}
        \exp\left\{
        iu^{\mathsf T}m-\frac12u^{\mathsf T}\Sigma u
        \right\}
        &=\E e^{iu^{\mathsf T}X_1} \\
        &=
        \E\left[
        \exp\left\{
        iu^{\mathsf T}M_\eta
        -\frac12u^{\mathsf T}(\Sigma-C)u
        \right\}
        \right].
\end{align*}
Consequently,
\[
        \E e^{iu^{\mathsf T}M_\eta}
        =
        \exp\left\{
        iu^{\mathsf T}m-\frac12u^{\mathsf T}Cu
        \right\},
        \qquad u\in\R^d.
\]
Thus \(M_\eta\sim\Normal(m,C)\).

Finally, by \eqref{eq:eta-vector-gaussian}, conditionally on \(M_\eta\), the
variables \(X_1,X_2,\ldots\) are i.i.d. with common law
\(\Normal(M_\eta,\Sigma-C)\), and \(M_\eta\sim\Normal(m,C)\).  Hence
\[
        (X_i)_{i\ge1}
        \stackrel{d}{=}
        (M+\xi_i)_{i\ge1},
\]
where
\[
        M\sim\Normal(m,C),\qquad
        \xi_1,\xi_2,\ldots \text{ are } \iid\ \Normal(0,\Sigma-C),
\]
and \(M\) is independent of \((\xi_i)_{i\ge1}\).  Therefore
\((X_i)_{i\ge1}\) is a Gaussian process.
\end{proof}

\begin{remark}
There is also a standard multidimensional form of Hardy's
uncertainty principle.  With the Fourier convention
\[
        \wh f(t)=\int_{\R^d}e^{it^{\mathsf T}x}f(x)\dd x,
\]
let \(A\) be a real symmetric positive definite matrix.  If \(f\in L^1(\R^d)\)
satisfies
\[
        f(x)
        =
        O\!\left(\exp\left\{-\frac12 x^{\mathsf T}A^{-1}x\right\}\right),
        \qquad
        \wh f(t)
        =
        O\!\left(\exp\left\{-\frac12 t^{\mathsf T}At\right\}\right),
\]
then
\[
        f(x)
        =
        \kappa\exp\left\{-\frac12 x^{\mathsf T}A^{-1}x\right\}
        \quad\text{a.e.}
\]
for some constant \(\kappa\in\mathbb C\).  This is the matrix version of the
same Gaussian borderline appearing in Theorem~\ref{thm:hardy}; see, for
example, Bonami--Demange--Jaming
\cite[Proposition~3.3]{BonamiDemangeJaming2003}. For the related multidimensional form and further background, see also Sitaram--Sundari--Thangavelu
\cite[Theorem~4]{SitaramSundariThangavelu1995} and
Fern\'andez-Bertolin--Vega~\cite{FernandezBertolinVega2025}.

Theorem~\ref{thm:vector-extension} could also be proved directly from a
multidimensional Hardy principle.  Namely, one would
localize the law \(\eta*\widetilde\eta\) on events
\(E\in\sigma(\eta)\), as in Section~\ref{sec:proof-difference-law}, and
apply Hardy's theorem to the resulting finite measures.  If \(2(\Sigma-C)\)
is positive definite, this gives
\[
        \eta*\widetilde\eta
        =
        \Normal(0,2(\Sigma-C))
        \qquad\text{a.s.}
\]
by the same argument as in the scalar proof.  If \(2(\Sigma-C)\)
is singular, the same route can be carried out after first reducing to the
support subspace of the corresponding Gaussian law.  (The previous projection proof avoids this separate degeneracy discussion.)
\end{remark}

\begin{funding}
\emergencystretch=2em
YZ is supported by the China Postdoctoral Science Foundation (Nos. 2023M743721 and 2025T180850).  QZ is supported by the National Natural Science Foundation of China (No. 12501162), the China Postdoctoral Science Foundation (No. 2024M762280), and the Jiangsu Provincial Natural Science Foundation (No. BK20250832).
\end{funding}

\end{document}